\documentclass[12pt,leqno]{amsart}
\usepackage{amssymb}
\usepackage{amscd}
\usepackage{amsmath}
\usepackage{amsthm}
\def\stacksum#1#2{{\stackrel{{\scriptstyle #1}}
{{\scriptstyle #2}}}}
\newcommand{\IGNORE}[1] {}
\newcommand{\eps}{\varepsilon}
\def\peter#1{\langle #1\rangle}
\def\ov#1{\overline{#1}}

\newcommand{\Q}{\mathbb Q}

\newcommand{\C}{\mathbb C}

\newcommand{\B}{\mathcal B}
\newcommand{\E}{\mathcal E}

\newcommand{\mfa}{\mathfrak{a}}

\newcommand{\mcO}{{\mathcal O}}
\newcommand{\Pic}{{\mathrm {Pic}}}
\newcommand{\Lalg}{L^{\mathrm{alg}}}
\newcommand{\Snew}{{\mathcal F}_{2k}(N)}
\newcommand{\Stnew}{{\mathcal F}_{2}(N)}
\newcommand{\modu}{\mathrm{mod}}

\makeatletter
\def\l@section{\@tocline{1}{4pt}{1pc}{}{}}
\def\l@subsection{\@tocline{2}{0pt}{2pc}{5pc}{}}
\makeatother

\newcommand\cyr{%
\renewcommand\rmdefault{wncyr}%
\renewcommand\sfdefault{wncyss}%
\renewcommand\encodingdefault{OT2}%
\normalfont \selectfont} \DeclareTextFontCommand{\textcyr}{\cyr}

\begin{document}
\newtheorem{Theorem}{Theorem}
\newtheorem*{TheoremA}{Theorem A}
\newtheorem*{corA}{Corollary A}
\newtheorem*{TheoremB}{Theorem B}
\newtheorem*{TheoremC}{Theorem C}
\newtheorem*{TheoremD}{Theorem D}
\newtheorem*{TheoremE}{Theorem E}
\newtheorem{Example}{Example}[section]
\newtheorem{Lemma}{Lemma}[section]
\newtheorem{Definition}{Definition}
\newtheorem{Proposition}{Proposition}
\newtheorem{Corollary}{Corollary}
\newtheorem{Conjecture}{Conjecture}
\newtheorem{Hypothesis}{Hypothesis}

\theoremstyle{definition}
\newtheorem*{note}{Note}
\newtheorem*{Remark}{Remark}
\title{Consequences of the Gross/Zagier formulae: \\ Stability of average $L$-values, subconvexity, and non-vanishing mod $p$}
\author{Philippe Michel and Dinakar Ramakrishnan}
\date{}

\maketitle

\medskip

\begin{flushright}
{\it In memory of Serge Lang}
\end{flushright}

\pagestyle{myheadings} \markboth{Philippe Michel and Dinakar
Ramakrishnan}{Consequences of the Gross/Zagier formulae}

\bigskip

\section*{Introduction}

\bigskip
In this paper we investigate some consequences of the Gross/Zagier
type formulae which were introduced by Gross and Zagier and then
generalized in various directions by Hatcher, Zhang, Kudla and
others \cite{GZ,Gross,Hatcher1,Zhang,MSRI}. Let us now recall these
formulae in the classical context. Denote by $K$ an imaginary
quadratic field of discriminant $-D$ say, with associated quadratic
character $\chi_{-D}=(\tfrac{-D}{\cdot})$, $\Psi$ a character of the
 ideal class group $\Pic(\mcO_{K})$ of $K$, ${\mathcal H}$ the upper
 half plane,
 and $g_{\Psi}$ the weight one theta series associated with $\Psi$:
 $$g_{\Psi}(z)=\sum_{m\geq 0}r_{\Psi}(m)q^m, \, \, q=\exp(2\pi \iota z), z\in {\mathcal H},$$
 where for $m\geq 1$
 $$r_{\Psi}(m)=\sum_{N(\mfa)=m}\Psi(\mfa)$$
 and $\mfa\subset\mcO_{K}$ ranging over the $\mcO_{K}$-ideal of norm
 $m$. We will denote the trivial character of $\Pic(\mcO_{K})$ by
 $1_K$.

 Now let $f$ be an holomorphic new cuspform of level $N$ coprime with  $D$,
  trivial nebentypus and weight $2k$:
$$f(z)=\sum_{m\geq 1}a_{m}(f)q^{m}.$$
Depending on how the primes dividing $N$ split in $K$, the
Gross/Zagier formula expresses the central value at $s=k$ (or the
derivative of that value) of the Rankin-Selberg $L$-function
$$L(s,f,\Psi):=L(2s,\chi_{-D})\sum_{m\geq 1}a_{m}(f)r_{\Psi}(m)m^{-s}$$
in term of an intersection/height pairing of the $f$-isotypic
component $e_{\Psi,f}$ of a cycle $e_{\Psi}$
 living on some Hecke module $M=M_{k,N}$: Denoting this pairing by
 $\peter{\cdot,\cdot}_{M}$ and the Petersson inner product on $S_{2k}(N)$ by
 $$\peter{f,g}=\int_{Y_{0}(N)}=f(z)\ov {g(z)}y^{2k-2}dxdy,$$ where $Y_0(N)$ denotes
 the open modular curve $\Gamma_0(N)\backslash {\mathcal H}$, one has
   \begin{equation}\label{GZformula}
   c_{k,K}\frac{L^{(i)}(k,f,\Psi)}{\peter{f,f}}=\peter{e_{\Psi,f},e_{\Psi,f}}_{M}
   \end{equation}
   for some constant $c_{k,K}>0$ and the order of derivative
   $i=i_{K,N}$ is $0$ or $1$ (depending on the sign of the functional equation).
   Originally
   the formula was proven as follows (for $i=0$):
let $M_{2k}(N)$ (resp. $S_{2k}(N)$) denote
   the space of holomorphic forms (resp. cusp forms) of weight $2k$ level $N$ and
   trivial nebentypus.
   The map $$f\mapsto L(s,f,\Psi)$$ being linear on $S_{2k}(N)$, can be represented
   by a kernel $f\mapsto \peter{f,G_{\Psi}}$ for some $G_{\Psi}\in M_{2k}(N)$ (same
   for the first derivative).
    By the Rankin-Selberg theory
    $$L(k,f,\Psi)=\int_{Y_{0}(N)}f(z)g_{\Psi}(z)E_{2k-1}(z)y^{(2k+1)/2-2}dxdy$$
    for a suitable holomorphic Eisenstein series
$E_{2k-1}$ of weight $2k-1$. The determination
    of $G_{\Psi}$ amounts to first taking the trace from level $N^\prime={\sl lcm}(4,N)$
    to $N$, and then computing the projection of $g_{\Psi}(z)E_{2k-1}(z)$ on
    $M_{2k}(N)$. This can be done and one infers from the
    computation of the Fourier expansion of $g_{\Psi}(z)E_{2k-1}(z)$, that the Fourier
    coefficients $a_{m}(G_{\Psi})$ of $G_{\Psi}$ are relatively elementary
    expressions involving the
    arithmetical functions $r_{\Psi}$ and variants thereof: see below for an example.
    One the other hand, using the theory of complex multiplication,
    Gross and Zagier, and subsequently other people, showed by an auxiliary computation that
    $$G_{\Psi}(z)=a_{0}(G_{\Psi})+\sum_{m\geq 1}\peter{T_{m}e_{\Psi},e_{\Psi}}_{M}q^m$$
    where $T_{m}$ denote the $m$-th Hecke operator acting on the module $M$. The final
    result follows then from a formal argument involving
     the multiplicity one theorem. The main observation underlying this paper is that
     the above computation provides formally an expression for the {\it average}
     of the central values $L(k,f,\Psi)$. Namely, if $S^{new}_{2k}(N)$ denote the set of
     arithmetically normalized new forms, then $\{f/\peter{f,f}^{1/2}\}_{f\in S^{new}_{2k}(N)}$
     may be completed to an orthonormal basis of $S_{2k}(N)$. Then decomposing $G_{\Psi}$ along
     such an orthonormal basis, and taking the $m$-th Fourier coefficient in the above
     decomposition,
     one deduces, for any $m\geq 1$,
     %\begin{multline*}\sum_{f\in S^{new}_{2k}(N)}\frac{L(k,f,\Psi)}{\peter{f,f}}a_{m}(f)
     %=a_{m}(G_{\Psi})\\
     %+{\mathcal Av}_{\rm old}+{\mathcal Av}_{\rm Eis}\\
     %\end{multline*}
     $$
     \sum_{f\in
     S^{new}_{2k}(N)}\frac{L(k,f,\Psi)}{\peter{f,f}}a_{m}(f) \, = \,
     a_{m}(G_{\Psi}) + +{\mathcal A}_{\rm old}(m)+{\mathcal A}_{\rm
     Eis}(m),
     $$
     where ${\mathcal A}_{\rm old}(m)$, resp. ${\mathcal A}_{\rm Eis}(m)$, is the contribution from the old forms,
     resp. the Eisenstein series, of weight $2k$ and level $N$.
     In principle, the $\hbox{Eisenstein series contribution}$ could be evaluated explicitly, while the
     $\hbox{ old forms contribution}$ could be computed by induction on $N$
     by  following the same scheme, though there is an added complication of finding a suitable orthonormal basis.
     We shall consider here the nicest possible situation for which
     these additional contributions
      have a particularly simple expression, in fact where the old part vanishes! Therefore we
      obtain,
     by the first step of the proof of the Gross/Zagier type formulae, a simple expression
     for the first moment
     $$\sum_{f\in S^{new}_{2k}(N)}\frac{L(k,f,\Psi)}{\peter{f,f}}a_{m}(f).$$
     Let us turn to a more specific example. Set $h=h_{K}=|\Pic(\mcO_{K})|$, the class number of
     $K$, $u=|\mcO_{K}^\times/\{\pm 1\}|$, and
     $$R(m):=\begin{cases}h/2u, \, &\mbox{ $m=0$}\\
     \sum\limits_\stacksum{\mfa\subset\mcO_{K}}{N(\mfa)=m}\, 1, \, &\hbox{$m\geq 1$}
     \end{cases},$$
     Moreover extend, for any ideal class group character $\Psi$, the definition of $r_{\Psi}(m)$
     to $m=0$ by setting
     $$r_{\Psi}(0)=\begin{cases}0, &\hbox{if $\Psi\not=1_K$}\\
     h/2u, &\hbox{if $\Psi=1_K$}.
     \end{cases}$$
     We also set
     $$\sigma_{N}(m)=\sum_\stacksum{d|m}{(d,N)=1}d$$
     Specializing to a generalization by Hatcher \cite{Hatcher1,Hatcher2} of a formula of Gross
     \cite{Gross}, we obtain
     \begin{Theorem}\label{identity} Let $-D<0$ be an odd fundamental discriminant; let $N$ be a
     prime which inert in $K=\Q(\sqrt{-D})$ and let $k\geq 2$ be an even integer.
      For $\Psi$ a character of $\Pic({\mcO}_K)$, then for
any positive integer $m$, we have the following exact identity:
\begin{multline*}
\noindent(2) \, \quad \quad
\frac{(2k-2)!D^{1/2}u^2}{2\pi(4\pi)^{2k-1}}\sum\limits_{f \in
{\mathcal F}_{2k}(N)}
\frac{L(f, \Psi, k)}{\peter{f,f}}a_m(f) \, = \,\\
-\delta\frac{12h^2}{N-1}{\sigma_N(m)} +{um^{k-1}r_\Psi(m)h} +
u^2m^{k-1}\sum_{n=1}^{\frac{mD}{N}}\Phi_k(n,\Psi,N)
\end{multline*}
Here
$$
\Phi_k(n,\Psi,N) \, = \, d((n,D))\delta_1(\Psi)R(n)r_\Psi(mD-nN)P_{k-1}(1-\frac{2nN}{mD}),
$$
with $P_{k-1}$ denoting the $(k-1)$-th Legendre polynomial;
$\delta\in \{0,1\}$ is $1$ iff $(k,\Psi)=(1,1_K)$;
$\delta_1(\Psi)\in \{0,1\}$ is $1$ if $D$ is prime, and when $D$ is
composite, it is $1$ iff $\Psi^2=1_K$ and there exist ideals $\frak
a, \frak b$, of respective norms $mD-nN$ and $n$, such that, for a
prime ideal $Q$ congruent to $-N$ mod $D$, the class of $\frak
a\frak bQ$ is a square in $\Pic(\mcO_{K})$.
\end{Theorem}
An asymptotic formula involving the average on the left was first
established for $k=1, \Psi=1_K$ by W.~Duke (\cite{Duke}), which
spurred a lot of other work, including that of Iwaniec and Sarnak
(\cite{IwaniecSarnak}) relating it to the problem of Siegel zeros
for $L(s,\chi_{-D})$. In the work of the second named author with
J.~Rogawski (\cite{RaRo}), a different proof of Duke's result was
given (for all weights), using Jacquet's relative trace formula
involving the integration of the kernel over the square of the split
torus, and in addition, the intervening measure was identified.

It is important to note that one obtains a {\it stability theorem}
when $N$ is sufficiently large compared with $D$ and $m$, and this
could perhaps be considered the most unexpected consequence of our
approach. Indeed, when $N>mD$, the sum on the far right of the
identity furnished by Theorem $1$ becomes zero, and our exact
average simplifies as follows:
\begin{Corollary} ({\rm Stability}) \, With the above notations and assumptions,
suppose moreover $N>mD$, then one has
\begin{multline*}
\frac{(2k-2)!D^{1/2}u^2}{2\pi(4\pi)^{2k-1}}\sum\limits_{f \in {\mathcal F}_{2k}(N)}
\frac{L(f, \Psi, k)}{\peter{f,f}}a_m(f) =\\
-\delta\frac{12h^2}{N-1}{\sigma_N(m)} +{um^{k-1}r_\Psi(m)h}
\end{multline*}
\end{Corollary}
We call the range $N>mD$, the {\it stable range}. As one can check
with other instances of the Gross/Zagier formulas, such as for the
derivative in the case of odd order of vanishing, this phenomenon
appears to be quite general. It has been recently generalized to
Hilbert modular forms of square-free level by B.~Feigon and
D.~Whitehouse (\cite{FW}), using the relative trace formula, now by
integrating the kernel over a non-split torus.

\medskip

When $\Psi=1_K$, we have the factorization
$$L(s,f,1_K)=L(s,f_K)=L(s,f)L(s,f\otimes \chi_{-D}),$$
where $f_K$ denotes the base change of $f$ to $K$, $L(s,f)$ the
Hecke $L$-function of $f$, and $f\otimes\chi_{-D}$ the twist of $f$
by $\chi_{-D}$. Thus for $m=1$ and $N>D$, we get the following
explicit identity involving the class number of $K$:
$$\frac{(2k-2)!D^{1/2}u}{2\pi(4\pi)^{2k-1}}\sum\limits_{f \in {\mathcal F}_{2k}(N)}
\frac{L(k,f)L(k,f\otimes\chi_{-D})}{\peter{f,f}}
={h}\bigl(1-\delta\frac{12h}{u(N-1)}\bigr)$$ In the weight 2 case,
as $N$ is taken to be a prime here, the cardinality of ${\mathcal
F}_2(N)$ is just the genus $g_0(N)$ of the compactification $X_0(N)$
of $Y_{0}(N)$. It is amusing to note that when $g_{0}(N)$ is zero,
one finds that
$$h=\frac{(N-1)u}{12},$$
implying that $h = 1$ when $(-D,N)$ is $(-3,5)$, $(-7,13)$,
$(-8,13)$ or $(-11,13)$, agreeing with known data. Similarly,
$X_0(11)$ is an elliptic curve $E/\Q$, and if we denote by $E_{-D}$
the $-D$-twist of $E$, we see, for $D=3$, that the algebraic special
value $A(1,E)A(1,E_{-3})$ is just $1/5$. In general one gets more
complicated identities, involving average central values, which are
all compatible with the Birch and Swinnerton-Dyer conjecture for
$E$, $E_{-D}$, and the Shafarevich-Tate groups {Sh}$(E)$,
{Sh}$(E_{-D})$.

\subsection{Application to the subconvexity problem}
We now discuss some simple applications of the above exact average
formula, the first one being a subconvex estimate for the central
values $L(k,f,\Psi)$. We refer to \cite{GAFA2000} for a general
discussion on the subconvexity problem. In the present case the
convexity bound is given by
$$
L(k,f,\Psi)\ll_{\eps}(kND)^\eps kN^{1/2}D^{1/2},\leqno(3)
$$
for any $\eps>0$. We prove here
\begin{Corollary} \, Preserve the notations of Theorem \ref{identity}.
Then for any $\eps>0$, we have
$$L(k,f,\Psi)\ll_{\eps}(kDN)^\eps kN^{1/2}D^{1/2}\bigl(\frac{1}{N^{1/2}}+\frac{N^{1/2}}{D^{1/2}}\bigr).$$
In particular this improves on convexity as long as
$$(kD)^\delta \leq N\leq D(kD)^{-\delta}$$
for some fixed $\delta>0$.
\end{Corollary}

Note that this breaks convexity for any fixed $k$, as long as $N$ is
between $D^\delta$ and $D^{1-\delta}$. The beauty is that we can
also vary $k$ in an appropriate region, obtaining a {\it hybrid
subconvexity}.

At this point we do not know of any application to these subconvex
estimates, but we are intrigued by them because they come for free
and seem to be hard to prove with the current methods of analytic
number theory (eg. see \cite{DFI,KMV}). Note also that such bounds
are fundamentally limited to the critical center $s=k$. For a
generalization to the Hilbert modular case, where $\Psi$ is allowed
to be any ray class character, see \cite{FW}.

\subsection{Application to non-vanishing problems}
Another line of application addresses the existence of $f$ for which
$L(k,f,\Psi)$ does not vanish. Indeed several variants of such
problems have been considered in the past by various methods
\cite{Duke,IwaniecSarnak,KM,OnoSkinner,Vatsal2}. Here we obtain
non-vanishing results which are valid with a fairly large uniformity
in the parameters, and again such uniformity seems hard to achieve
by purely analytic methods.
\begin{Theorem}Assumptions being as for Theorem A. Suppose that $$N\gg_{\delta} D^{1/2+\delta}$$
for some $\delta>0$, then there exists $f\in S^{new}_{2k}(N)$ such that
$$L(k,f,\Psi)\not=0.$$
The same conclusion holds as long as $N>D$ and either $k\not=1$ or
$\Psi\not=1_K$.
\end{Theorem}
When $\Psi=1_K$, we also obtain non-vanishing result in a somewhat
greater range:
\begin{Theorem} Suppose $\Psi=1_K$, $k=1$ and
$$h<\frac{N-1}{12}.$$ Then there exist $f$ such that
$$L(k,f)L(k,f\otimes\chi_{-D})\not=0.$$
\end{Theorem}
  Non-vanishing theorems of this kind, with an {\em explicit} dependence between $N$ and $D$
  (like $N>D$ or $N-1>12h$), are of some interest. For instance, in the paper
 \cite{Merel1}, Merel needs to consider the following problem: Given a prime $p$
 and a character $\chi$ of conductor $p$ which is not even and quadratic, does there exist an
 $f\in\mathcal{F}_{2}(p)$ such that $L(1,f\otimes\chi)\not=0$? In the appendix of that paper, the first
 named author and E. Kowalski prove that this is the case when $p$ is greater than an
  explicit but very large number. In particular, it has so far not been possible to answer the problem numerically
  in the finitely many remaining cases; this has been answered however for
   $p<1000$ \cite{Merel2}.
Closer to the main concern of the present paper, Ellenberg
\cite{Ellenberg1,Ellenberg2}
  uses analytic methods to prove the non-vanishing of the twisted $L$-function
  $L(1,f\otimes\chi_{-4})$ for some $f$ in $\Stnew$
  for $N$ of the form $p^2$ or $2p^2$ ($p$ an odd prime) and with prescribed eigenvalues at the
  Atkin/Lehner operators $w_{2},w_{p}$,
  subject to an {\em explicit} lower bound on $p$. Ellenberg concludes from this the
  non-existence of primitive integral solutions to the
   generalized Fermat equation
   $A^4+B^2=C^{p}$
   as long as $p>211$; that this equation has only a finite number of primitive solutions is a theorem of
   Darmon and Granville. Another related set of examples is in the work of Dieulefait and Urroz (\cite{DU}).
   In a sequel to this paper under preparation (\cite{MiRa}), we will develop a suitable generalization of the
   exact average formula to a class of composite levels $N$, and investigate similar questions by modifying the method.
   This extension is subtle for three reasons:  $N$ is not square-free,
   $D$ is not odd, and $N,D$ are not relatively prime.

     \subsection{Nonvanishing modulo $p$}
     The exactness of the Gross/Zagier formulae even enable us to obtain {\it average non-vanishing
     results} for the {\it algebraic part} of the $L(k,f,\Psi)$ modulo suitable primes $p$.
      Again, such a question has been considered in the past, see for example
      \cite{BJOSK,Vatsal2}. However, these earlier works
      addressed the question of the existence of the non-vanishing of
      $L(k,f,\Psi)$ mod $p$  when the form $f$ is {\em fixed} and when the character $\Psi$ varies. Here our results go
      in the other direction as we fix $p$ and let $N$ and $f$ vary. Given
      $f\in\Snew$ and $g_{\Psi}$ as above, we denote by $\Lalg(k,f,\Psi)$ the algebraic part of
      $L(k,f,\Psi)$ (see section 5, (11), for a precise definition). It follows from the work of Shimura
      that $\Lalg(k,f,\Psi)$ is an algebraic number satisfying the reciprocity law
      $$\Lalg(k,f,\Psi)^\sigma=\Lalg(k,f^\sigma,\Psi^\sigma)$$
      for any $\sigma$ automorphism of $\C$ \cite{Shimura}.
      \begin{Theorem}\label{padic} Let $p>2k+1$ be a prime, $\mathcal P$ be a chosen place in $\ov{\Q}$
      above $p$ and let $N,D$ be as in Theorem \ref{identity}.
       Suppose moreover that $p$ does not divide $h=h_{-D}$, that $N>D$, and that $N$ is greater that some
       absolute constant. Then there exists $f\in \Snew$ such that
      $$\Lalg(k,f,\Psi)\not\equiv 0 \, \, (\modu\ {\mathcal P}).$$
      \end{Theorem}
Naturally, the question of integrality of $\Lalg(k,f,\Psi)$, which
is subtle, and our result only concerns the numerator of the
$L$-value. When $\Psi=1_K$, we also prove the following variant:

             \begin{Theorem}\label{padic2}Notations and assumptions as in Theorem \ref{padic}. Suppose
             moreover that $\Psi=1$ and $N>pD$.
              Then there exists $f\in \Snew$ such that
      $$\sqrt{D}(2\pi)^{-2k}\frac{L(k,f)L(k,f\otimes\chi_{-D})}{\langle f,f\rangle}a_{p}(f)\not\equiv 0 \, \, (\modu\ {\mathcal P}^{2k-1}).$$
      \end{Theorem}
      The assertion makes sense because the left hand side is (see
      section 5.1) a $p$-unit times
      $a_p(f)$ times $\Lalg(k,f,1_K)$.

      \medskip

There are two fundamental periods $c^+(f)$ and $c^-(f)$ associated
to $f$ such that for any Dirichlet character $\nu$, the special
value $\Lalg(k,f\otimes\nu)$, defined as $L(k,f\otimes \nu)/c^{{\rm
sgn}(\nu(-1))}(f)$ times a simple factor (see section 5, (12)) is an
algebraic number. One gets the near-factorization
$$
\eta_f\Lalg(k,f,1_K) \, = \, \Lalg(k,f)\Lalg(k,f\otimes\chi_{-D}),
$$
where $\eta_f$ is essentially the order of the congruence module
considered by Hida, Wiles, Taylor, Flach, Diamond, and others, which
measures the congruences $f$ has with other modular forms modulo
$p$. The needed non-divisibility properties of $\eta_f$ (for
suitable $p$) are understood (at least) if $f$ is ordinary or $k=1$.
Now finally, let us suppose we are in the classical weight $2$
situation, i.e., with $\Psi=1_K$ and $k=1$.
\begin{Theorem}\label{padic3} Let $p$ an odd prime not dividing $Dh_{-D}$,
with $D$ odd. Then there exist infinitely many newforms of $f$ of
prime level $N$ and weight $2$ such that
$$
{\rm num}\left(\frac{\Lalg(1,f\otimes\chi_{-D} )}{\eta_f}\right) \,
\not\equiv \, 0 \, \pmod p,
$$
where $\eta_f$ is the order of the congruence module of $f$.
      \end{Theorem}

See section 5 for a discussion of $\eta_f$, which measures the
congruences which $f$ may have with other modular forms of the same
weight and level. An analogue of Theorem 6 should also hold, in a
suitable range of $p$, for forms of higher weight, and this question
will be taken up elsewhere.

\medskip

\subsection{Acknowledgement} Serge Lang always conveyed infectious excitement about
Mathematics to anyone he came into contact with, and he will be
missed. He was quite interested in the values of $L$-functions and
in the {\it divisibility properties} of arithmetic invariants, and
it is a pleasure to dedicate this article to him. The first author
would like to thank Caltech for its hospitality during the
preparation of this work. The second author would like to thank
Flach, Hida, Prasanna and Vatsal for helpful conversations
concerning the last part, and the National Science Foundation for
support through the grant DMS0402044.

\section{The weight $2$ case}

It may be instructive to explain why the exact average formula holds
in the weight $2$ case when $\Psi=1$. Let $B$ be a quaternion
division algebra over $\Q$, ramified only at $N$ and $\infty$, with
maximal order $R$. Put $Y$ is the associated rational curve such
that Aut$(Y) = B^\ast/\Q^\ast$. Put
$$
X = B^\ast\backslash Y \times \hat{B}^\ast/\hat{R}^\ast =
\cup_{j=1}^n \Gamma_j\backslash Y,
$$
where $\hat{B}^\ast=\prod\limits_p{}' B_p^\ast$ and
$\hat{R}^\ast=\prod\limits_p R_p^\ast$, with each $\Gamma_j$ being a
finite group. Then Pic$(X)$ identifies with $\{e_1, e_2, \ldots,
e_n\}$, where each $e_j$ is the class of $\Gamma_j\backslash Y$.
{Since} $N$ is inert in $K=\Q[\sqrt{-D}]$, there is an embedding
$f\in {\rm Hom}(K,B) = Y(K)$. It results in certain {\it Heegner
points} $x=(f,b)$ of discriminant $-D$ in $X$, with $b \in
\hat{B}^\ast/\hat{R}^\ast$. For any eigenform $f$, let $c_f$ denote
the $f$-component of $c = \sum_A x_A$, where $A$ runs over ideal
classes of $K$. Then by a beautiful theorem of B.~Gross ([G]),
providing an analogue for the $L$-value of the Gross-Zagier theorem
for the first derivative, one has
$$
\langle c_f, c_f\rangle \, = \,
u^2\sqrt{D}\frac{L(1,f)L(1,f\otimes\chi_{-D})}{(f,f)},
$$
where $\langle \cdot, \cdot\rangle$ is a natural {\it height
pairing} on Pic$(x)$. We have by orthogonality,
$$
\langle c,T_mc\rangle = \langle c_E,T_mc_E\rangle +\sum\limits_f
\langle c_f,T_mc_f\rangle,
$$
where $T_m$ is the operator corresponding to the $m$-the Hecke
operator on $M_2(N)$, $f$ runs over newforms in $M_2(N)$, and $E$
denotes the unique (holomorphic) Eisenstein series (of weight $2$
and level $N$). Using the fact that $f$ and $E$ are Hecke
eigenforms, and that $\langle c_E, c_E\rangle \, = \,
\frac{12h^2}{N-1}$, we get by averaging Gross's formula,
$$
u^2\sqrt{\vert D\vert}\sum\limits_f
\frac{L(1,f)L(1,f\otimes{\chi_{-D}})}{(f,f)} =
-\sigma_N(m)\frac{12h^2}{N-1} + \langle c, T_mc\rangle.
$$
One has
$$
\langle c,T_mc\rangle \, = \, \sum\limits_A\sum\limits_B \langle
x_B,T_mx_{AB}\rangle. $$ If we pick $q \equiv -N ($mod $D)$, with
$q\mcO_K = Q\overline Q$ in $K$, one sees that
$$
\sum\limits_B\langle x_B,T_mx_{AB}\rangle \, = \, uhR_A(m)
+\sum\limits_{n=1}^{mD/N} R_A(mD-nN)d((n,D))R_{\{QA\}}(n),
$$
with
$$
R_{\{QA\}}(n) = \vert\{I : N(I)=n, QAI \in {\rm
Pic}({\mcO}_K)^2\}\vert.
$$
Note that $R_{\{QA\}}(n)$ is just $R_A(n)$ when $D$ is prime. The
assertion of Theorem 1 now follows by summing over $A$. Moreover,
when $mD$ is less than $N$, $\sum\limits_B\langle
x_B,T_mx_{AB}\rangle$ simply equals $uhR_A(m)$, and this furnishes
Corollary 1 (stability) in the weight $2$ case.

\section{\bf Proof of the main identity for all $k\geq 1$}

\subsection{Preliminaries} \, For $N\geq 1$, let $M_{2k}(N)$ (resp
$S_{2k}(N)$) denote, as usual, the space of holomorphic modular
forms (resp. cusp forms) of weight $2k$ level $N$ and trivial
character. For $f\in M_{2k}(N)$, we write the Fourier expansion at
the infinite cusp as
\[
f(z)=\sum_{m\geq 0}a_m(f)q^m, q=\exp (2\pi\imath z).
\]
We denote by $\Snew$, the set of cuspidal new forms $f$ (normalized
in the usual way, so that the first Fourier coefficient $a_1(f)$ is
1. Whenever it converges, we denote the Petersson inner product on
$M_{2k}(N)$ by
\[
\langle f,g\rangle =\int_{Y_{0}(N)}f(z)\overline{g(z)}y^{2k-2}dxdy.
\]
Let $-D<0$ be an odd fundamental discriminant, $K=\Q (\sqrt{-D}), \mcO_k$ the
maximal order of $K , \Pic(\mcO_K)$ the ideal class group, and $u=u_k=|
\mcO_K{}^\times |/2$.  For any ideal class $A\in{\rm Pic}(\mcO_k)$, define
\[
r_A(m)=\begin{cases}|\{\mfa\subset\mcO_K,N(\mfa )=m,\mfa\in A\}| &\text{if }
m\geq 1\\
\frac{1}{2u} & \text{if }m=0
\end{cases}
\]
The theta series
\[
\theta_A(z)=\sum_{m\geq 0}r_A(m)q^m,q=\exp (2\pi\imath z)
\]
is a modular form of weight 1, level $D$ and central character
$\chi_{-D}$. Moreover, for any $\Psi\in\widehat{{\rm Pic}(\mcO_K)}$,
put
\[
\theta_\Psi (z)=\sum_A\overline\Psi (A)\theta_A(z),
\]
whose Fourier coefficients are then given by
$$
a_m(\theta_\Psi) = \sum_A \overline\Psi(A)a_m(\theta_A).
$$
In particular, the constant term $a_0(\theta_\Psi)$ equals
$\frac{1}{2u}\sum_A \overline\Psi(A)$, which is, by orthogonality,
zero iff $\Psi\ne 1_K$, when $\theta_\Psi$ is a cusp form. Setting
\[
L(s,f,A):=\sum_{\stacksum{n>1}{(n,N)=1}}\frac{\chi_{-D}(n)}{n^{1+1(s-k)}}\sum_{m\geq
1} \frac{a_m(f)r_A(m)}{m^s},
\]
one has
\[
L(s,f,\Psi )=\sum_{A\in\Pic (\mcO_K)}\Psi (A)L(s,f,A).
\]
Define a holomorphic function $G_A$ on the upper half plane
${\mathcal H}$, invariant under $z\rightarrow z+1$, by means of its
Fourier expansion at infinity:
\begin{equation}
G_A(z):=\sum^\infty_{m=0}b_{m,A}q^m,
\end{equation}
where
\begin{align}
b_{m,A}&=m^{k-1}\frac{h}{u}r_A(Dm)\\
&+m^{k-1}\sum^{mD/N}_{n=1}\delta (n)r_A(mD-nN)R_{(-NA)}(n)P_{k-1}\left (
1-\frac{2nN}{mD}\right )\nonumber
\end{align}
In this definition, $u$ and $R(n)=\sum_Ar_A(n)$ are as in the Introduction,
$\delta (n)$ is 1 (resp. 2) if $(m,D)$ is 1 (resp. $\neq 1$), and for $r\geq 0,
P_r$ is the $r$-th Legendre polynomial defined by
\[
P_r(x):=\frac{1}{2^r}\sum^{[r/2]}_{m=1}(-1)^m\begin{pmatrix}r\\m\end{pmatrix}
\begin{pmatrix}2r-2m\\r\end{pmatrix}x^{r-2m}.
\]

The following result, due to B. Gross, D. Zagier and R. Hatcher, is
crucial to us:
\begin{Theorem}
$G_A$ is a modular form of weight $2k$, level $N$, and trivial character;
it is cuspidal if $k>1$, and for every newform $f$ of weight $2k$ and level
$N$, we have
\[
L(k,f,A)=\frac{(4\pi )^{2k}}{2(2k-2)!D^{1/2}}(f,G_A).
\]
\end{Theorem}

For $k=1$, see [11], Prop. 9.1, and for general $k$, this is in [12],
Theorem 5.6 and [14], Theorem 5.2.  (See also [13], where the case $D$ prime
is treated.)

\subsection{The exact average formula} Let
\[
E \, = \, E_{2,N} \, = \, \sum^\infty_{n=0}a_n(E)q^n
\]
denote a holomorphic Eisenstein series for $\Gamma_0(N)$ of weight
2.  Since $N$ is prime, the modular curve $Y_{0}(N)$ has only two
cusps, namely $\infty$ and 0.  It then follows that $E$ is unique up
to scalar multiple, and so $E(z)/a_0(E)$ is well defined with
constant term 1 at $\infty$.  To be specific, we will take
\[
E(z)=\frac{N-1}{12}+\sum^\infty_{m=1}\sigma_N(m)q^m,
\]
where $\sigma_N(m)=\sum_{d|m,(d,N)=1}d$.

For $A\in\Pic (\mcO_K)$, with $G_A$ being as in the previous section, put
\begin{equation}
G^{\rm cusp}_A(z):G_A(z)-\delta_{k=1}\frac{b_{0,A}}{a_0(E)}E(z),
\end{equation}
with $\delta_{k,1}$ being 1 (resp. 0) if $k=1$ (resp. $k\neq 1$).  Then
$G^{\rm cusp}_A$ is a holomorphic cusp form of level $N$, weight $2k$,
and trivial character, with coefficients $a_m(G^{\rm cusp}_A)$.
\medskip

\noindent{\bf Lemma 2.1.} {\it For $-D$ an odd fundamental
discriminant and $N$ a prime inert in $K$, we have, for any $m\geq
1$,
\begin{align*}
\frac{2(2k-2)!D^{1/2}}{(4\pi )^{2k}}\sum_{f\in\Snew}&\frac{L(k,f,A)}{\langle
f,f\rangle}a_m(f)\\
&=a_m(G^{\rm cusp}_A)=b_{m,A}-\delta_{k=1}\frac{b_{0,A}}{a_0(E)}a_m(E)
\end{align*}}

In order to prove this, we first need the following
\medskip

\noindent{\bf Lemma 2.2.}  {\it Assume that $N$ is a prime which is
inert in $K=\Q [\sqrt{-D}]$.  Let $\varphi$ be any old form in
$S_{2k}(N)$.  Then we have, for every $A\in\Pic (\mcO_K)$,
\[
(\varphi ,G^{\rm cusp}_A)=0.
\]}

There is nothing to prove when $k<6$, since $S_{2k}(1)$ is zero in
that case (cf. \cite{L}, for example.)  Such a Lemma will not in
general hold for composite $N$.
\medskip

{\bf Proof of Lemma 2.2.}  Since $\varphi$ is cuspidal, it suffices
to prove that $(\varphi ,G_A)=0$.  Put
\[
G_\Psi :=\sum_{A\in\Pic (\mcO_K)}\Psi (A)G_A
\]
which is modular form of weight 1 and character $\chi_{-D}$.  It is sufficient
to show that $(\varphi ,G_\Psi )=0$ for all ideal class characters $\Psi$ of $K$.
If $\varphi =\sum^\infty_{n=1}a_n(\varphi )q^n$, put
\begin{equation}
D(s,\varphi \times\theta_\Psi )=\sum^\infty_{n=1}\frac{a_n(\varphi )\overline a_n
(\theta_\Psi)}{n^s}
\end{equation}
Then the Rankin-Selberg method give the identity
\begin{equation}
(2\pi )^{-k}\Gamma (k)D(k,\varphi\times\theta_\Phi )=\langle f,{\rm Tr}_{ND/N}
(\theta_\Phi\E_{2k-1,N})\rangle
\end{equation}
where $\E_{2k-1,N}$ is the result of slashing a holomorphic
Eisentein series of weight $2k-1$ (and character $\chi_{-D}$) with
the Atkin involution $u_N$, and Tr$_{ND/D}$ denotes the trace from
$S_{2k}(ND)$ to $S_{2k}(N)$.  In fact, the calculations of Gross and
Zagier (\cite{GZ}) show that
\begin{equation}
G_\Psi ={\rm Tr}_{ND/N}(\theta_\Psi\E_{2k-1,N}).
\end{equation}
Now let $\varphi$ be a newform of level 1 (and weight $2k$).  Then
since $N$ is prime, it defines two old forms of level $N$, namely
$\varphi_1(z)=\varphi (z)$ and $\varphi_2(z)=\varphi (Nz)$, so that
$a_m(\varphi_2)$ is zero unless $N|m$, and
$a_mN(\varphi_2)=a_m(\varphi )$.  Since the new and old forms are
orthogonal to each other under $(\cdot ,\cdot )$, and since the
space of old forms of level $N$ are spanned by $\{\varphi_d,d=1,N\}$
with $\varphi$ running overl all the cusp forms of level 1, it
suffices to prove that each $D(k, \varphi_d\times\theta_\Psi )=0$.
Let $d=1$.  Then one obtains (by section 3, Lemma 1, of [Sh]):
\begin{equation}
L(2k,\chi_{-D})D(k,\varphi_d\times\theta_\Psi )=L(k,\varphi\times\theta_\Psi ).
\end{equation}
Since $L(x,\chi_{-D})$ is non-zero at $s=2k$ (which is in the region of
absolute convergence), it reduces to checking the vanishing of the right hand
side.  Since $\varphi$ has level 1, the root number of $L(k,\varphi \times
\theta_\Psi )$ is $-1$, yielding the requisite vanishing.  When $d=N,D(k,
\varphi_d\times\theta_\Psi )$ is still a non-zero multiple of $L(k,\varphi\times
\theta_\Psi )$, which is zero.

\hfill$\Box$
\medskip

{\bf Proof of Lemma 2.1}  We may choose an orthogonal basis $\B$ of
$S_{2k}(N)$ to be of the form $\Snew\cup\B'$, where $\B'$ consists
of old forms.  Clearly we have
\begin{equation}
\sum_{f\in\B}\frac{(f,G^{\rm cusp}_A)}{\langle f,f\rangle}=G^{\rm cusp}_A.
\end{equation}
In view of the Lemma, the sum on the left hand side needs to run
only over newforms $f$.  Applying Theorem 6, and using (8), we
obtain
\[
\frac{2(2k-2)!D^{1/2}}{(4\pi )^{2k}}\sum_{f\in\Snew}\frac{L(f,\Psi ,k)}
{\langle f,f\rangle}=G^{\rm cusp}_A.
\]
The lemma now follows by taking the $m$-coefficient of the above identity.

\hfill$\Box$
\medskip

{\bf Proof of Theorem 1}  The exact average formula follows by
performing the averaging $\sum_{A\in\Pic (\mcO_K)}\Psi (A)\dots$ on
both sides of the formula in Lemma 2.1 using the formula (5) for the
coefficients $b_{m.A}$, and by noting that
\[
\frac{a_m(E)}{a_0(E)}=\frac{12}{N-1}\sigma_N(m)
\]
and that $b_{0,A}=\tfrac{h}{2u^2}$.

\hfill$\Box$

\section{\bf Subconvex Bounds}

In this section, we prove Corollary 2.  By the work of Waldspurger,
Guo and Jacquet (\cite{Guo, Waldspurger}; also \cite{KohZ} for
$\Psi=1_K$),
\[
L(k,f,\Psi )\geq 0.
\]
Thus from formula (2) for $m=1$, we have
\[
\frac{(2k-2)!D^{1/2}}{2(4\pi )^{2k}}\frac{L(f,\Psi ,k)}{\langle f,f\rangle}
\leq\frac{h}{u}+\sum^{\frac{D}{N}}_{n=1}|\Psi_k(n,\Psi ,N)|
\]
Since $|P_{k-1}(x)|\leq 1$ for $|x|\leq 1$ and $R(n),|r_\Psi (n)|\leq d(n)$
(where $d(n)$ is the number of divisors of $n$), so that
\[
R(n)|r_\Psi (D-nN)|\leq d(n)^2+d(D-nN)^2,
\]
we see that the $n$-sum on the right side is bounded by $\tfrac{D}{N}(\log
D)^3$.  From the class number formula, we have
\[
h\ll D^{1/2}\log D
\]
and
\[
\langle f,f\rangle\ll (4\pi )^{-2k}(2k-1)!N(\log kN)^3
\]
as follows from \cite{ILS}, (2.3), (unlike the corresponding bound
for Maass forms (\cite{HL}) this upper bound is elementary since $f$
holomorphic so its Fourier coefficients satisfy the
Ramanujan|Petersson bound). Thus we see that
\[
L(f,\Psi ,k)\ll (\log kN)^3(\log D)^3k(N+D^{1/2}).
\]
\hfill$\Box$

\section{\bf Application to non-vanishing}

We prove here Theorem 2.  Arguing exactly as above we have
\begin{align*}
\frac{(2k-2)!D^{1/2}}{2(4\pi )^{2k}}\sum_{f\in\Snew}\frac{L(f,\Psi
,k)}{\langle f,f
\rangle}&=\frac{h}{u}-\delta\frac{6(h/u)^2}{N-1}+O\left (\frac{D}
{N}(\log D)^3\right )\\
&=\frac{h}{u}+O\left (\frac{D}{N}(\log D)^3\right )
\end{align*}
By Siegel's Theorem, which gives $h=D^{1/2 +o(1)}$, we see that the
right side is positive as soon as $N>D^{1/2+\delta}$ for some
$\delta >0$. If $N>D$, then we are in the stable range and we have
\begin{equation}
\frac{(2k-2)!D^{1/2}}{2(4\pi )^{2k}}\sum_{f\in\Snew}\frac{L(f,\Psi ,k)}
{\langle f,f\rangle}=\frac{h}{u}\left (1-\delta\frac{6(h/u)}{N-1}
\right ).
\end{equation}
When $\delta =0$, this concludes the proof of Theorem 2 since $h\geq
1$. \hfill$\Box$
\medskip

Suppose now that $\delta =1$ (ie. $k=1,\Psi =1_K$). Then we remark
that
\[
\sum^{\frac{D}{N}}_{n=1}\Psi_1(n,1,N)\geq 0
\]
so that
\[
\frac{(2k-2)!D^{1/2}}{2(4\pi )^{2k}}\sum_{f\in\Snew}\frac{L(f,\Psi ,k)}
{\langle f,f\rangle}\geq\frac{h}{u}\left (1-\frac{6(h/u)}{N-1}
\right )
\]
combining the proof of Theorem 3.\hfill$\Box$

\section{\bf Non-vanishing mod $p$}

\subsection{\bf Algebraic Parts of $L$-values}

Let us put
\begin{equation}
\Lalg(k,f, \Psi) \, = \,
(-1)^k(2\pi)^{-2k}(k-1)!^2g(\chi_{\-D})\frac{L(k,f,\Psi)}{\langle f,
f\rangle},
\end{equation}
where $g(\overline{\Psi})$ is the Gauss
sum. Then it is known, by Shimura (\cite{Shimura}, see also
\cite{Hd1}), that $\Lalg(k,f, \psi)$ is an algebraic number obeying
the reciprocity law:
$$
\Lalg(k,f^\sigma,\Psi^\sigma) = \Lalg(k,f, \Psi)^\sigma,
$$
for every automorphism $\sigma$ of $\C$.

Next recall that for $\Psi=1_K$, $L(k,f,\Psi)$ factors as
$L(k,f)L(k,f\otimes\chi_{-D})$. For any Dirichlet character $\nu$,
the algebraic part of $L(k,f\otimes\nu)$ is given by
\begin{equation}
\Lalg(k,f\otimes\nu) \, = \, g(\overline
\nu)(k-1)!\frac{L(k,f,\nu)}{(-2\pi i)^kc_\pm(f)},
\end{equation}
where $c_\pm(f)$ is a fundamental period of $f$, with $\pm =
\nu(-1)$. Again, one has for any automorphism $\sigma$ of $\C$,
$\Lalg(k,f^\sigma\otimes\nu^\sigma)$ is
$\Lalg(k,f\otimes\nu)^\sigma$.

This leads to the near-factorization
\begin{equation}
\eta_f\Lalg(k,f,1_K) \, = \, \Lalg(k,f)\Lalg(k,f\otimes\chi_{-D}),
\end{equation}
where $\eta_f$ equals, thanks to a series of papers of
Hida (cf. \cite{Hd1}, \cite{Hd2}), Wiles (\cite{Wiles}),
Taylor-Wiles (\cite{TW}), and Diamond-Flach-Guo (\cite{DFG}), the
order of the congruence module of $f$, i.,e the number which counts
the congruences of $f$ with other modular forms of the same weight
and level.

\subsection{\bf Proof of Theorems 4 and 5}
From the definition of the algebraic part, the hypothesis of Theorem
4 and the formula (9), used in conjunction with $\delta =0$, we have
(up to multiplication by a $p$-unit)
\[
\sum_{f\in\Snew}L^{\rm alg}(k,f,\Psi )=\frac{h}{u}.
\]
The conclusion of Theorem 4 is immediate.

For the proof of Theorem 5, we have, assuming that $N>pD$,
\[
\sum_{f\in\Snew}L^{\rm alg}(k,f,1_K)=\frac{h}{u}\left
(1-\frac{12(h/u)}{N-1} \right ).
\]
Therefore the conclusion holds except possibly if
$p|(1-\tfrac{6(h/u)}{N-1})$. Suppose we are in that latter case.
Then we apply the exact formula of Corollary 1 with $m=p$ and get
\[
\sum_{f\in\Snew}L^{\rm alg}(k,f,1_K)a_p(f)=\frac{h}{u}\left (R(p)-
\frac{6(h/u)}{N-1}(p+1)\right )
\]
$R(p)$ is either 0 or 2, if it is zero, then the left hand side of
the previous formula is not divisible by $p$. If $R(p)=2$, then
$2-\tfrac{6(h/u)}{N-1}$ is not divisible by $p$ since by assumption
$p|(1-\tfrac{6(h/u)}{N-1})$. So we are done in all
cases.\hfill$\Box$

\medskip

\subsection{Proof of Theorem 6}

\medskip

Here are restricting to the weight $2$ case, and by the theory of
modular symbols, cf. Stevens \cite{St} and Vatsal \cite{V} - see
also Prasanna \cite{P} - we know that for any Dirichlet character
$\nu$, the special value $\Lalg(1,f\otimes\nu)$ is integral except
possibly at the Eisenstein primes; these are the primes dividing
$$
\tilde{N}: = \, \prod_{q\vert N} q(q^2-1),
$$
which is related to the order of the cuspidal divisor class group,
studied for modular curves, among others, by Kubert and Lang.

We may, and we will, choose $N$ to lie in the infinite family of
primes which are inert in $K$ and are such that $p \nmid \tilde{N}$.

Now Theorem 6 follows by the near-factorization (13) of
$\Lalg(1,f,1_K)$. It may be useful to note that when $f$ has
$\Q$-coefficients, with associated elliptic curve $E$ over $\Q$, one
knows (cf. Flach \cite{F}) that any prime dividing $\eta_F$ also
divides the degree of the modular parametrization $X_0(N) \to E$.

\vskip 0.2in

\bibliographystyle{math}    % replacing \bib...{plain}
\bibliography{Mi-Ram}

\end{document}